\documentclass[11pt]{article}
\usepackage{amssymb}
\def\o{\'o}
\def\n{\'n}
\def\z{{\accent95 z}}
\def\a{\kern+.6ex\lower.42ex\hbox{$\scriptstyle \iota$}\kern-1.20ex a}

\newcommand{\bbK}{{\mathbb K}}

\newcommand{\bm}{{\mathbf m}}
\newcommand{\bx}{{\mathbf x}}

\newcommand{\g}{\mathbf g}
\newcommand{\h}{\mathbf h}

\def\rank{\mbox{\rm rank}}
\def\ord1{\mbox{\rm ord}}
\def\mod{\mbox{\rm mod}}
\def\P{\partial}

\title{Formal and convergent solutions of analytic equations
\footnotetext{
\begin{minipage}[t]{5.2in}
{\small
MSC: 14B12\\
Keywords: formal solution, parametric solution
} 
\end{minipage}
}}
\begin{document}
\author{by Arkadiusz P\l{}oski}
\date{}
\maketitle
\begin{abstract}\noindent
{\small
We provide the detailed proof of a strengthened version of the M.~Artin
Approximation Theorem.}
\end{abstract}
\begin{center}
\begin{minipage}{4.2in}\small
Impressed by the power of the Preparation Theorem --
indeed, it prepares us so well! -- I considered
``Weierstrass Preparation Theorem and its immediate
consequences'' as a possible title for the entire book.

\vspace{2ex}\hspace{7em}\noindent Sheeram S.~Abhyankar, Preface to
\cite{Abh1964}
\end{minipage}
\end{center}

\section{Introduction}

The famous Approximation Theorem of M.~Artin~\cite{Art1968}
asserts that any formal solution of a system of analytic equations
can be approximated by convergent solutions up to a given order.
In my PhD thesis~\cite{Pl1973} I was able by analysis of the argument
used in~\cite{Art1968} to sharpen the Approximation Theorem: any formal
solution can be obtained by specializing parameters in a convergent
parametric solution. The theorem was announced with a sketch of proof
in~\cite{Pl1974}. The aim of theses notes is to present the
detailed proof of this result. It is based on the Weierstrass
Preparation Theorem. The other tools are: a Jacobian Lemma which
is an elementary version of the Regularity Jacobian Criterion used in
~\cite{Art1968}, the trick of Kronecker (introducing and specializing
variables) and a generalization of the Implicit Function Theorem due
to Bourbaki~\cite{Bou1961} and Tougeron~\cite{Tou1972}.
All theses ingredients are vital in the proofs of some other results
of this type (see~\cite{Art1969},~\cite{Wav1975}).
For more information on approximation theorems in local analytic geometry
we refer the reader to Teissier's article~\cite{Te1993-94} and to
Chapter~8 of the book~\cite{deJong-Pf2000}.

Let $\bbK$ be a field of characteristic zero with a non-trivial valuation.
We put $\bbK[[x]]=\bbK[[x_1,\dots,x_n]]$ the ring of formal power series
in variables $x=(x_1,\dots,x_n)$ with coefficients in $\bbK$.
If $f=\sum_{k\geq p}f_k$ is a nonzero power series represented as the sum
of homogeneous forms with $f_p\neq 0$ then we write $\ord1\,f=p$.
Additionally we put $\ord1\,0=+\infty$ and use the usual conventions
on the symbol $+\infty$. The constant term of any series $f\in\bbK[[x]]$
we denote by $f(0)$. A power series $u\in\bbK[[x]]$ is a unit if $uv=1$
for a power series $v\in\bbK[[x]]$. Note that $u$ is a unit if and only if
$u(0)\neq 0$. The non-units of $\bbK[[x]]$ form the unique maximal ideal
$\bm_{\bx}$ of the ring $\bbK[[x]]$. The ideal $\bm_{\bx}$ is generated
by the variables $x_1,\dots,x_n$. One has $f\in\bm_{\bx}^c$, where $c>0$
is an integer, if and only if $\ord1\,f\geq c$. Recall that if
$g_1,\dots,g_n\in\bbK[[y]]$, $y=(y_1,\dots,y_n)$ are without constant
term then the series $f(g_1,\dots,g_n)\in\bbK[[y]]$ is well-defined.
The mapping which associates with $f\in\bbK[[x]]$ the power series
$f(g_1,\dots,g_n)$ is the unique homomorphism sending $x_i$ for
$g_i$ for $i=1,\dots,n$. Let $\bbK\{x\}$ be the subring of
$\bbK[[x]]$ of all convergent power series. Then $\bbK\{x\}$ is
a local ring. If $g_1,\dots,g_n\in\bbK\{y\}$ then
$f(g_1,\dots,g_n)\in\bbK\{y\}$ for any $f\in\bbK\{x\}$.

In what follows we use intensively the Weierstrass Preparation and
Division Theorems. The reader will find the basic facts concerning
the rings of formal and convergent power series in \cite{Abh1964},
\cite{Lef1953} and \cite{Zar-Sam1960}.

Let $f(x,y)=(f_1(x,y),\dots,f_m(x,y))\in\bbK\{x,y\}^m$ be convergent power
series in the variables $x=x_1,\dots,x_n)$ and $y=(y_1,\dots,y_N)$ where
$m,n,N$ are arbitrary non-negative integers.
The theorem quoted below is the main result of \cite{Art1968}.

\vspace{1ex}\noindent{THE ARTIN APPROXIMATION THEOREM}\\
{\it Suppose that there exists a sequence of formal power series
$\bar y(x)=(\bar y_1(x),\dots,\bar y_N(x))$ without constant term
such that
$$
  f(x,\bar y(x))=0\;.
$$
Then for any integer $c>0$ there exists a sequence of convergent power
series $y(x)=(y_1(x),\dots,y_N(x))$ such that
$$
  f(x,y(x))=0\mbox{ and }y(x)\equiv\bar y(x)\;(\mod\,\bm_{\bx}^c)\;.
$$
} 

\vspace{1ex}\noindent The congruence condition means that the power
series $y_\nu(x)-\bar y_\nu(x)$ are of order$\mbox{}\geq c$ i.e. the
coefficients of monomials of degree$\mbox{}<c$ agree in
$y_\nu(x)$ and $\bar y_\nu(x)$. We will deduce the Artin Approximation
Theorem from the following result stated with a sketch of proof
in \cite{Pl1974}.

\vspace{1ex}\noindent{THEOREM}
{\it With the notation and assumptions of the Artin theorem there
exists a sequence of convergent power series
$y(x,t)=(y_1(x,t),\dots,y_N(x,t))\in\bbK\{x,t\}^N$, $y(0,0)=0$, where
$t=(t_1,\dots,t_s)$ are new variables, $s\geq 0$, and a sequence of
formal power series
$\bar t(x)=(\bar t_1(x),\dots,\bar t_s(x))\in\bbK[[x]]^s$, $\bar t(0)=0$
such that
$$
  f(x,y(x,t))=0\mbox{ and }\bar y(x)=y(x,\bar t(x))\;.
$$
} 

\vspace{1ex}\noindent The construction of the parametric solution $y(x,t)$
depends on the given formal solution $\bar y(x)$. To get the Artin
Approximation Theorem from the stated above result fix an integer $c>0$.
Let $y(x,t)$ and $\bar t(x)$ be series such as in the theorem and let
$t(x)=(t_1(x),\dots,t_s(x))\in\bbK\{x\}^s$ be convergent power series
such that $t(x)\equiv\bar t(x)\;\mod\,\bm_{\bx}^c$. Therefore
$y(x,t(x))\equiv y(x,\bar t(x))\;\mod\,\bm_{\bx}^c$ and it suffices to
set $y(x)=y(x,t(x))$~\rule{1ex}{1ex}

Before beginning the proof of the theorem let us indicate two corollaries
of it.

\vspace{1ex}\noindent{COROLLARY 1}.
{\it Assume that $m=N$, $f(x,\bar y(x))=0$ and
$$
\det\frac{J(f_1,\dots,f_N)}{J(y_1,\dots,y_N)}(x,\bar y(x))\neq 0\;.
$$
Then the power series $\bar y(x)$ are convergent.
} 

\vspace{1ex}\noindent Proof. Let $y(x,t)$ and $\bar t(x)$ be power series
without constant term such that $f(x,y(x,t))=0$ and
$\bar y(x)=y(x,\bar t(x))$. It is easy to check by differentiation
of equalities $f(x,y(x,t))=0$ that $(\P{y_\nu}/\P{t_\sigma})(x,t)=0$
for $\nu=1,\dots,N$ and $\sigma=1,\dots,s$. Therefore the series
$y(x,t)$ are independent of $t$ and the series $\bar y(x)$ are
convergent~\rule{1ex}{1ex}

\vspace{1ex}\noindent{COROLLARY 2}.
{\it If $f(x,y)\in\bbK\{x,y\}$ is a nonzero power series of $n+1$
variables $(x,y)=(x_1,\dots,x_n,y)$ and $\bar y(x)$ is a formal power series
without constant term such that $f(x,\bar y(x))=0$ then $\bar y(x)$ is
a convergent power series.
} 

\vspace{1ex}\noindent Proof. By Corollary~1 it suffices to check
that there exists a power series $g(x,y)\in\bbK\{x,y\}$ such that
$g(x,\bar y(x))=0$ and $(\P{g}/\P{y})(x,\bar y(x))\neq 0$.
Let $I=\{g(x,y)\in\bbK\{x,y\}:\,g(x,\bar y(x))=0\}$. Then
$I\neq\bbK\{x,y\}$ is a prime ideal of $\bbK\{x,y\}$. Assume the contrary,
that is, that for every $g\in I$: $(\P{g}/\P{y})\in I$. Then we get
by differentiating the equality $g(x,\bar y(x))=0$ that
$(\P{g}/\P{x_i})\in I$ for $i=1,\dots,n$ and, by induction, all
partial derivatives of $g$ lie in $I$. Consequently $g=0$ for
every $g\in I$ i.e. $I=(0)$. A contradiction since
$0\neq f\in I$~\rule{1ex}{1ex}

\section{Reduction to the case of simple solutions}

We keep the notation introduced in Introduction. We will call a sequence
of formal power series $\bar y(x)\in\bbK[[x]]$, $\bar y(0)=0$ a
{\it simple solution\/} of the system of analytic equations $f(x,y)=0$ if
$f(x,\bar y(x))=0$ and
$$
\rank\frac{J(f_1,\dots,f_m)}{J(y_1,\dots,y_N)}(x,\bar y(x))=m\;.
$$
Thus, in this case, $m\leq N$.

In what follows we need

\vspace{1ex}\noindent{THE JACOBIAN LEMMA}.
{\it Let $I$ be a nonzero prime ideal of the ring $\bbK\{x\}$,
$x=(x_1,\dots,x_n)$. Then there exist covergent power series
$h_1,\dots,h_r\in I$ such that
\begin{itemize}
\item[\rm(i)] $\displaystyle
\rank\frac{J(h_1,\dots,h_r)}{J(x_1,\dots,x_n)}(\mod\,I)=r\;,$
\item[\rm(ii)] $\forall h\in I$, $\exists a\notin I$ such that
$ah\in(h_1,\dots,h_r)\bbK\{x\}$.
\end{itemize}
} 

\vspace{1ex}\noindent Before proving the above lemma let us note that it is
invariant with respect to $\bbK$-linear nonsingular transformations.
If $\Phi$ is an authomorphism of $\bbK\{x\}$ defined by
$$
\Phi(f(x_1,\dots,x_n))=
f\left(\sum_{j=1}^nc_{1j}x_j,\dots,\sum_{j=1}^nc_{nj}x_j\right)
$$
with $\det(c_{ij})\neq 0$ then the Jacobian Lemma is true for $I$ if and
only if it is true for $\Phi(I)$.

\vspace{1ex}\noindent Proof of the Jacobian Lemma (by induction
on the number $n$ of variables $x_i$).\\
If $n=1$ then $I=(x_1)\bbK\{x_1\}$ and $h_1=x_1$.
Suppose that $n>1$ and that the lemma is true for prime ideals
of the ring of power series in $n-1$ variables.
Using a $\bbK$-linear nonsingular transformation we may assume
that the ideal $I$ contains a power series $x_n$-regular of order $k>0$
i.e. such that the term $x_n^k$ appears in the power series with a
non-zero coefficients. Therefore, by the Weierstrass Preparation Theorem
$I$ contains a distinguished polynomial
$$
  w(x',x_n)=x_n^k+a_1(x')x_n^{k-1}+\dots+a_k(x'),\mbox{ where }
  x'=(x_1,\dots,x_{n-1})\;.
$$
By the Weierstrass Division Theorem every power series $h=h(x)$ is of the
form $h(x)=q(x)w(x',x_n)+r(x',x_n)$ where $r(x',x_n)$ is an
$x_n$-polynomial (of degree$\mbox{}<k$). Therefore, the ideal $I$
is generated by the power series which are polynomials in $x_n$
and to prove the Jacobian Lemma it suffices to find power series
$h_1,\dots,h_r$ such that (i) holds and (ii) is satisfied for
$h\in I\cap\bbK\{x'\}[x_n]$.

Let $I'=I\cap\bbK\{x'\}$ and consider the set $I\setminus I'[x_n]$.
Clearly $w(x',x_n)\in I\setminus I'[x_n]$. Let
$$
  h_1(x',x_n)=c_0(x')x_n^l+c_1(x')x_n^{l-1}+\dots+c_l(x')
$$
be a polynomial in $x_n$ of the minimal degree $l$, $l\geq 0$,
which belongs to $I\setminus I'[x_n]$. Since the degree $l\geq 0$
is minimal, we have
$$
\begin{array}{ccc}
l & > & 0\;,\\
c_0(x') & \notin & I'\;,\\
\frac{\P{h_1}}{\P{x_n}} & \notin & I'\;.\\
\end{array}
$$
Let $h(x',x_n)\in I$ be a polynomial in $x_n$. Dividing $h(x',x_n)$
by $h_1(x',x_n)$ (Euklid's division) we get
\begin{itemize}
\item[(E)] $c_0(x')^p\,h(x',x_n)=q(x',x_n)\,h_1(x',x_n)+r_1(x',x_n)$
\end{itemize}
where $x_n$-degree of $r_1(x',x_n)$ is less than $l$ and $p\geq 0$
is an integer. Since the $x_n$-degree of $r_1(x',x_n)$ is$\mbox{}<l$
then all coefficients of $r_1(x',x_n)$ lie in $I'$. If $I'=(0)$
then $r_1(x',x_n)=0$ and (E) proves the Jacobian Lemma.

If $I'\neq (0)$ then by the induction hypothesis there exists series
$h_2,\dots,h_r\in I'$ such that
\begin{itemize}
\item[(i$'$)] $\displaystyle
\rank\frac{J(h_2,\dots,h_r)}{J(x_1,\dots,x_{n-1})}(\mod\,I')=r-1\;,$
\item[(ii$'$)] $\forall h'\in I'$, $\exists a'\notin I'$ such that
$a'h'\in(h_2,\dots,h_r)\bbK\{x'\}$.
\end{itemize}
We claim that $h_1,\dots,h_r$ satisfy (i) and (ii) of the Jacobian Lemma.
To check (i) observe that
$$
\det\frac{J(h_1,\dots,h_r)}{J(x_{i_1},\dots,x_{i_{r-1}},x_n)}
=\det\frac{J(h_2,\dots,h_r)}{J(x_{i_1},\dots,x_{i_{r-1}})}\cdot
\frac{\P{h_1}}{\P{x_n}}\;,
$$
where $i_1,\dots,i_{r-1}\in\{1,\dots,n-1\}$ and use (i$'$).
Applying (ii$'$) to the coefficients of $r_1(x',x_n)$ we find a power
series $a'(x')$ such that $a'(x')\,r_1(x',x_n)\in (h_2,\dots,h_r)\bbK\{x\}$.
By (E) we get $a(x')\,h(x',x_n)\in(h_1,\dots,h_r)\bbK\{x\}$ where
$a(x')=a'(x')\,c_0(x')^p\notin I$ which proves (ii)~\rule{1ex}{1ex}

\vspace{1ex}\noindent Now, we can check

\vspace{1ex}\noindent{PROPOSITION 2.1.}
{\it Let $f(x,y)=(f_1(x,y),\dots,f_m(x,y))\in\bbK\{x,y\}^m$, $f(x,y)\neq 0$,
$\bar y(x)=(\bar y_1(x),\dots,\bar y_N(x))\in\bbK[[x]]$, $\bar y(0)=0$,
be formal power series such that $f(x,\bar y(x))=0$. Then there exist
convergent power series $h(x,y)=(h_1(x,y),\dots,h_r(x,y))\in\bbK\{x,y\}^r$
such that
\begin{itemize}
\item[\rm(i)] $h(x,\bar y(x))=0$,
\item[\rm(ii)] $\displaystyle
\rank\frac{J(h_1,\dots,h_r)}{J(y_1,\dots,y_{N})}(x,\bar y(x))=r\;,$
\item[\rm(iii)] suppose that there exist power series
$y(x,t)=(y_1(x,t),\dots,y_N(x,t))$, $y(0,0)=0$ and
$\bar t(x)=(\bar t_1(x),\dots,\bar t_N(x))$, $\bar t(0)=0$, such that
$h(x,y(x,t))=0$ and $\bar y(x)=y(x,\bar t(x))$. Then $f(x,\bar y(x))=0$.
\end{itemize}
} 

\vspace{1ex}\noindent Proof. Consider the prime ideal
$$
I=\{g(x,y)\in\bbK\{x,y\}:\,g(x,\bar y(x))=0\}\;.
$$
Clearly $f_1(x,y),\dots,f_m(x,y)\in I$ and $I\neq(0)$.
By the Jacobian Lemma there exist series $h_1(x,y),\dots,h_r(x,y)\in I$
such that
\begin{itemize}
\item $\displaystyle
\rank\frac{J(h_1,\dots\dots\dots,h_r)}{J(x_1,\dots,x_n,y_1,\dots,y_N)}
(x,\bar y(x))=r\;,$
\item $\forall g\in I$, $\exists a\notin I$ such that
$a(x,y)\,g(x,y)\in(h_1,\dots,h_r)\bbK\{x,y\}$.
\end{itemize}
We claim that $h_1,\dots,h_r$ satisfy the conditions (i), (ii), (iii).
Condition (i) holds since $h_1,\dots,h_r\in I$. To check (ii) it suffices
to observe that
\begin{itemize}
\item[(J)] $\displaystyle
\rank\frac{J(h_1,\dots\dots\dots,h_r)}{J(x_1,\dots,x_n,y_1,\dots,y_N)}
(x,\bar y(x))=
\rank\frac{J(h_1,\dots,h_r)}{J(y_1,\dots,y_N)}(x,\bar y(x))\;.$
\end{itemize}
Indeed, differentiating the equations $h_i(x,\bar y(x))=0$, $i=1,\dots,r$,
we get
$$
\frac{\P{h_i}}{\P{x_j}}(x,\bar y(x))+\sum_{\nu=1}^N
\frac{\P{h_i}}{\P{y_\nu}}(x,\bar y(x))\frac{\P{\bar y_\nu}}{\P{x_j}}=0
\mbox{ for }j=1,\dots,n
$$
and (J) follows. To check (iii) let us write
$$
a_i(x,y)f_i(x,y)=\sum_{k=1}^ra_{i,k}(x,y)h_k(x,y)\mbox{ in }\bbK\{x,y\}\;,
$$
where $a_i(x,y)\notin I$ for $i=1,\dots,m$. Thus $a_i(x,\bar y(x))\neq 0$
and $a_i(x,\bar y(x,t))\neq 0$ since $\bar y(x)=y(x,\bar t(x))$ and (iii)
follows~\rule{1ex}{1ex}

\section{The Bourbaki--Tougeron implicit function theorem}

Let $f(x,y)=(f_1(x,y),\dots,f_m(x,y))\in\bbK\{x,y\}^m$ be convergent
power series in the variables $x=(x_1,\dots,x_n)$ and $y=(y_1,\dots,y_N)$.
Suppose that $m\leq N$ and put
$$
J(x,y)=\frac{J(f_1,\dots,f_m)}{J(y_{N-m+1},\dots,y_N)}\mbox{ and }
\delta(x,y)=\det\,J(x,y)\;.
$$
Let $M(x,y)$ be the adjoint of the matrix $J(x,y)$. Thus we have
$$
M(x,y)J(x,y)=J(x,y)M(x,y)=\delta(x,y)I_m\;
$$
where $I_m$ is the identity matrix of $m$ rows and $m$ columns.
Let $g(x,y)=(g_1(x,y),\dots,g_m(x,y))\in\bbK\{x,y\}^m$ be convergent
power series defined by
$$
\left[\begin{array}{c}g_1(x,y)\\\vdots\\g_m(x,y)\end{array}\right]=
M(x,y)\left[\begin{array}{c}f_1(x,y)\\\vdots\\f_m(x,y)\end{array}\right]\;.
$$
It is easy to see that
\begin{itemize}
\item[(a)] $g_i(x,y)\in(f_1(x,y),\dots,f_m(x,y))\bbK\{x,y\}$ for $i=1,\dots,m$
\end{itemize}
and
\begin{itemize}
\item[(b)] $\delta(x,y)f_i(x,y)\in(g_1(x,y),\dots,g_m(x,y))\bbK\{x,y\}$ for $i=1,\dots,m$.
\end{itemize}
Now, we can state

\vspace{1ex}\noindent{THE BOURBAKI-TOUGERON IMPLICIT FUNCTION THEOREM}\\
{\it Suppose that there exists a sequence of formal power series
$y^0(x)=(y_1^0(x),\dots,y_N^0(x))$, $y^0(0)=0$, such that
$$
g(x,y^0(x))\equiv 0\;\mod\,\delta(x,y^0(x))^2\bm_{\bx}\;.
$$
Then
\begin{itemize}
\item[\rm I.] Let $y_\nu(x,t)=y^0_\nu(x)+\delta(x,y^0(x))^2t_\nu$ for
$\nu=1,\dots,N-m$ where $t=(t_1,\dots,t_{N-m})$ are new variables.
Then there exists a unique sequence of formal power series
$u(x,t)=(u_{N-m+1}(x,t),\dots,u_N(x,t))\in\bbK[[x,t]]^m$, $u(0,0)=0$,
such that if we let $y_\nu(x,t)=y_\nu^0(x)+\delta(x,y^0(x))u_\nu(x,t)$
for $\nu=N-m+1,\dots,N$ and $y(x,t)=(y_1(x,t),\dots,y_N(x,t))$ then
$$
  f(x,y(x,t))=0\mbox{ in }\bbK[[x,t]]\;.
$$
If the series $y^0(x)$ are convergent then $u(x,t)$ and $y(x,t)$
are covergent as well.

\item[\rm II.] For every sequence of formal power series
$\bar y(x)=(\bar y_1(x),\dots,\bar y_N(x))$, $\bar y(0)=0$, the following
two conditions are equivalent
\begin{itemize}
\item[\rm(i)] there exists a sequence of formal power series
$\bar t(x)=(\bar t_1(x),\dots,\bar t_{N-m}(x))$, $\bar t(0)=0$, such that
$\bar y(x)=y(x,\bar t(x))$,
\item[\rm(ii)] $f(x,\bar y(x))=0$ and
\begin{eqnarray*}
\bar y_\nu(x) & \equiv & y^0_\nu(x)\;\mod\,\delta(x,y^0(x))^2\bm_{\bx}^2
\mbox{ for }\nu=1,\dots,N-m\\
\bar y_\nu(x) & \equiv & y^0_\nu(x)\;\mod\,\delta(x,y^0(x))\bm_{\bx}^2
\mbox{ for }\nu=N-m+1,\dots,N\;.
\end{eqnarray*}
\end{itemize}
\end{itemize}
} 

\vspace{1ex}\noindent{Remark.} In what follows we call
\begin{center}
\begin{tabular}{ll}
$y^0(x)$ & an approximate solution of the system $f(x,y)=0$, \\
$y(x,t)$ & a parametric solution determined by the approximate solution $y^0(x)$ \\
$\bar y(x)$ & satifying (i) or (ii) a subordinate solution to the approximate solution $y^0(x)$ \\
\end{tabular}
\end{center}
Proof. Let $v=(v_1,\dots,v_N)$ and $h=(h_1,\dots,h_n)$ be variables.
Taylor's formula reads
\begin{itemize}
\item[(T)]
\begin{minipage}{5.5in}
\begin{eqnarray*}
\left[\begin{array}{c}f_1(x,v+h)\\\vdots\\f_m(x,v+h)\end{array}\right] & = &
\left[\begin{array}{c}f_1(x,v)\\\vdots\\f_m(x,v)\end{array}\right]+
\frac{J(f_1,\dots,f_m)}{J(y_1,\dots,y_{N-m})}(x,v)
\left[\begin{array}{c}h_1\\\vdots\\h_{N-m}\end{array}\right]\\
& + & J(x,v)\left[\begin{array}{c}h_{N-m+1}\\\vdots\\h_{N}\end{array}\right]+
\left[\begin{array}{c}P_{1}(u,v,h)\\\vdots\\P_{m}(u,v,h)\end{array}\right]
\end{eqnarray*}
\end{minipage}
\end{itemize}
where $P_i(x,v,h)\in(h_1,\dots,h_{N})^2\bbK\{x,v,h\}$ for $i=1,\dots,m$.
Let $u=(u_{N-m+1},\dots,u_N)$ be variables and put
{\small
$$
F_i(x,t,u)=f_i(x,y_1(x,t),\dots,y_{N-m}(x,t),
y^0_{N-m+1}(x)+\delta(x,y^0(x))u_{N-m+1},\dots,
y^0_{N}(x)+\delta(x,y^0(x))u_{N})\;.
$$
} 
Substituting in Taylor's formula (T) $v_i=y^0_i(x)$ for $i=1,\dots,N$,
$h_i=\delta(x,y^0(x))^2t_i$ for $i=1,\dots,N-m$ and
$h_i=\delta(x,y^0(x))u_i$ for $i=N-m+1,\dots,N$ we get
\begin{eqnarray*}
\left[\begin{array}{c}F_1(x,t,u)\\\vdots\\F_m(x,t,u)\end{array}\right] & = &
\left[\begin{array}{c}f_1(x,y^0(x))\\\vdots\\f_m(x,y^0(x))\end{array}\right]+
\delta(x,y^0(x))^2\frac{J(f_1,\dots,f_m)}{J(y_1,\dots,y_{N-m})}(x,y^0(x))
\left[\begin{array}{c}t_1\\\vdots\\t_{N-m}\end{array}\right]\\
& + & \delta(x,y^0(x))J(x,y^0(x))\left[\begin{array}{c}u_{N-m+1}\\\vdots\\u_{N}\end{array}\right]+
\delta(x,y^0(x))^2\left[\begin{array}{c}Q_{1}(x,t,u)\\\vdots\\Q_{m}(x,t,u)\end{array}\right]
\end{eqnarray*}
where $Q_i(x,t,u)\in(t,u)^2\bbK\{x,t,u\}$ for $i=1,\dots,m$.
Multiplying the above identity by the matrix $M(x,y^0(x))$ and taking
into account that $M(x,y^0(x))J(x,y^0(x))=\delta(x,y^0(x))I_m$ and
$g_i(x,y^0(x))\equiv 0$ $(\mod\,\delta(x,y_0(x))^2\bm_{\bx})$ for
$i=1,\dots,m$, we get
\begin{itemize}
\item[(*)] $\displaystyle
M(x,y^0(x))\left[\begin{array}{c}F_1(x,t,u)\\\vdots\\F_m(x,t,u)\end{array}\right]=
\delta(x,y^0(x))^2\left[\begin{array}{c}G_1(x,t,u)\\\vdots\\G_m(x,t,u)\end{array}\right]
$
\end{itemize}
where $G_i(0,0,0)=0$ for $i=1,\dots,m$. Differentiating (*) we obtain
$$
M(x,y^0(x))\frac{J(F_1,\dots,F_m)}{J(u_{N-m+1},\dots,u_N)}(x,t,u)=
\delta(x,y^0(x))^2\frac{J(G_1,\dots,G_m)}{J(u_{N-m+1},\dots,u_N)}(x,t,u)
$$
which implies
\begin{itemize}
\item[(**)] $\displaystyle
\det\frac{J(F_1,\dots,F_m)}{J(u_{N-m+1},\dots,u_N)}(x,t,u)=
\delta(x,y^0(x))^{m+1}det\frac{J(G_1,\dots,G_m)}{J(u_{N-m+1},\dots,u_N)}(x,t,u)
$
\end{itemize}
since $\det\,M(x,y^0(x))=\delta(x,y^0(x))^{m-1}$. On the other hand
$$
\frac{J(F_1,\dots,F_m)}{J(u_{N-m+1},\dots,u_N)}(x,0,0)=
\delta(x,y^0(x))^{m}J(x,y^0(x))
$$
and
$$
\det\frac{J(F_1,\dots,F_m)}{J(u_{N-m+1},\dots,u_N)}(x,0,0)=
\delta(x,y^0(x))^{m+1}\;.
$$
Therefore we get from (**)
$$
\det\frac{J(G_1,\dots,G_m)}{J(u_{N-m+1},\dots,u_N)}(x,0,0)=1\;,
$$
in particular
$$
\det\frac{J(G_1,\dots,G_m)}{J(u_{N-m+1},\dots,u_N)}(0,0,0)=1\;.
$$
By the Implicit Function Theorem there exist formal power series
$$
u(x,t)=(u_{N-m+1}(x,t),\dots,u_{N}(x,t))
$$
such that
{\small
$$
(G_1(x,t,u),\dots,G_m(x,t,u))\bbK[[x,t,u]]=
(u_{N-m+1}-u_{N-m+1}(x,t),\dots,u_{N}-u_{N}(x,t))\bbK[[x,t,u]]\;.
$$
} 
If $y^0(x)$ are convergent then $G(x,t,u)$ and $u(x,t)$ are convergent as well.
In particular $G(x,t,u(x,t))=0$ and by (*) $F(x,t,u(x,t))=0$ which implies
$f(x,y(x,t))=0$ where $y_\nu(x,t)=y^0_\nu(x)+\delta(x,y^0(x))u_\nu(x,t)$
for $\nu=N-m+1,\dots,N$. Let
$\tilde u(x,t)=(\tilde u_{N-m+1}(x,t),\dots,\tilde u_{N}(x,t))$,
$\tilde u(0,0)=0$, be power series such that $f(x,\tilde y(x,t))=0$
where
{\footnotesize
$$
\tilde y(x,t)=(y_1(x,t),\dots,y_{N-m}(x,t),y^0_{N-m+1}(x)+\delta(x,y^0(x))%
\tilde u_{N-m+1}(x,t),\dots,y^0_{N}(x)+\delta(x,y^0(x))\tilde u_{N}(x,t))\;.
$$
} 
Then $F(x,t,\tilde u(x,t))=0$ and by (*) $G(x,t,\tilde u(x,t))=0$. Thus
we get $\tilde u(x,t)=u(x,t)$. This proves the first part of the
Bourbaki-Tougeron Implicit Function Theorem. To check the second part it
suffices to observe that for any formal power series
$\bar t(x)=(\bar t_1(x),\dots,\bar t_{N-m}(x))$ and
$\bar u(x)=(\bar u_{N-m+1}(x),\dots,\bar u_{N}(x))$ without constant term
$G(x,\bar t(x),\bar u(x))=0$ if and only if
$\bar u(x)=u(x,\bar t(x))$~\rule{1ex}{1ex}

\section{Approximate solutions}

We keep the notions and assumptions of Section~3.

\vspace{1ex}\noindent{PROPOSITION 4.1}.
{\it Let $\bar y(x)=(\bar y_1(x),\dots,\bar y_N(x))$, $\bar y(0)=0$,
be a formal solution of the system of analytic equations $f(x,y)=0$,
$m\leq N$, such that the power series $\delta(x,\bar y(x))$ is
$x_n$-regular of strictly positive order $p>0$. Then there exists
an approximate solution $\bar v(x)\in\bbK\{x'\}[x_n]^N$ of the system
$f(x,y)=0$ and such that $\bar y(x)$ is a solution of $f(x,y)=0$
subordinate to $\bar v(x)$.
} 

\vspace{1ex}\noindent Proof. 
By the Weierstrass Preparation Theorem
$\delta(x,\bar y(x))=\bar a(x)\cdot\mbox{unit}$ where
$$
\bar a(x)=x_n^p+\sum_{j=1}^p\bar a_j(x')x_n^{p-j}
$$
is a distinguished polynomial. Using the Weierstrass Division Theorem we get
$$
\bar y_\nu(x)=\sum_{j=0}^{2p-1}\bar v_{\nu,j}(x')x_n^j+
\bar a(x)^2(c_\nu+\bar t_\nu(x))\mbox{ for }\nu=1,\dots,N-m
$$
and
$$
\bar y_\nu(x)=\sum_{j=0}^{p-1}\bar v_{\nu,j}(x')x_n^j+
\bar a(x)(c_\nu+\bar u_\nu(x))\mbox{ for }\nu=N-m+1,\dots,N
$$
where $c_\nu\in\bbK$ for $\nu=1,\dots,N$, while
$$
\bar t(x)=(\bar t_1(x),\dots,\bar t_{N-m}(x))\mbox{ and }
\bar u(x)=(\bar u_{N-m+1}(x),\dots,\bar u_{N}(x))
$$
are formal power series without constant term. Let
$$
\bar v_\nu(x)=\sum_{j=0}^{2p-1}\bar v_{\nu,j}(x')x_n^j+
\bar a(x)^2 c_\nu\mbox{ for }\nu=1,\dots,N-m
$$
and
$$
\bar v_\nu(x)=\sum_{j=0}^{p-1}\bar v_{\nu,j}(x')x_n^j+
\bar a(x)c_\nu\mbox{ for }\nu=N-m+1,\dots,N\;.
$$
Clearly $\bar v(x)=(\bar v_1(x),\dots,\bar v_N(x))\in\bbK[[x']][x_n]^N\;.$

\vspace{1ex}\noindent{Propery 1:}
{\it $\delta(x,\bar v(x))=\bar a(x)\cdot\mbox{unit}$ } 

\vspace{1ex}\noindent Proof. From $\bar y(x)\equiv\bar v(x)$
$(\mod\,\bar a(x)\bm_{\bx})$ we get
$\delta(x,\bar y(x))\equiv\delta(x,\bar v(x))$ $(\mod\,\bar a(x)\bm_{\bx})$
and Property~1 follows since
$\delta(x,\bar y(x))=\bar a(x)\cdot\mbox{unit}$~\rule{1ex}{1ex}

\vspace{1ex}\noindent{Propery 2:}
{\it $g(x,\bar v(x))\equiv 0$ $(\mod\,\bar a(x)^2\bm_{\bx})$} 

\vspace{1ex}\noindent Proof. Substituting in Taylor's formula (T)
$v=\bar v(x)$, $h_\nu=\bar a(x)^2\bar t_\nu(x)$ for $\nu=1,\dots,N-m$
and $h_\nu=\bar a(x)\bar u_\nu(x)$ for $\nu=N-m+1,\dots,N$ we get
\begin{eqnarray*}
\left[\begin{array}{c}0\\\vdots\\0\end{array}\right] & = &
\left[\begin{array}{c}f_1(x,\bar v(x))\\\vdots\\f_m(x,\bar v(x))\end{array}\right]+
\bar a(x)^2\frac{J(f_1,\dots,f_m)}{J(y_1,\dots,y_{N-m})}(x,\bar v(x))
\left[\begin{array}{c}\bar t_1(x)\\\vdots\\\bar t_{N-m}(x)\end{array}\right]\\
& + & \bar a(x)J(x,\bar v(x))\left[\begin{array}{c}\bar u_{N-m+1}(x)\\\vdots\\\bar u_{N}(x)\end{array}\right]+
\bar a(x)^2\left[\begin{array}{c}\bar Q_{1}(x)\\\vdots\\\bar Q_{m}(x)\end{array}\right]\;.
\end{eqnarray*}
Multiplying the above identity by $M(x,\bar v(x))$ and taking into account
the formula
$$M(x,\bar v(x))J(x,\bar v(x))=\delta(x,\bar v(x))I_m$$
we get Property~2~\rule{1ex}{1ex}

\vspace{1ex}\noindent{PROPOSITION~4.2}
{\it Let $(c^0_{\nu,j})$, $\nu=1,\dots,N$, $j=0,1,\dots,D$, be a family
of constants such that $c^0_{\nu,0}=0$ for $\nu=1,\dots,N$. Suppose that
$$
\left(\sum_{j=0}^Dc_{1,j}^0x_n^j,\dots,\sum_{j=0}^Dc_{N,j}^0x_n^j\right)
$$
is an approximate solution of the system of equations $f(0,x_n,y)=0$
such that
$$
\ord1\,\delta\left(0,x_n,\sum_{j=0}^Dc_{1,j}^0x_n^j,\dots,\sum_{j=0}^Dc_{N,j}^0x_n^j\right)
=p\,,\quad 0<p<+\infty\;.
$$
Let $V^0=(V^0_{\nu,j})$, $\nu=1,\dots,N$, $j=0,1,\dots,D$ be variables.
Then there exists a sequence
$$
F(x',V^0)=(F_1(x',V^0),\dots,F_M(x',V^0))\in\bbK\{x',V^0\}^M
$$
such that for any family $(\bar v^0_{\nu,j}(x'))$ of formal power series
without constant term the following two conditions are equivalent
\begin{itemize}
\item[\rm(i)] $\displaystyle
\left(\sum_{j=0}^D(c_{1,j}^0+\bar v^0_{1,j}(x'))x_n^j,\dots,
\sum_{j=0}^D(c_{N,j}^0+\bar v^0_{N,j}(x'))x_n^j\right)$\\
is an approximate solution of the system $f(x,y)=0$,
\item[\rm(ii)] $F(x',(\bar v^0_{\nu,j}(x')))=0$ in $\bbK[[x']]$.
\end{itemize}
} 

\vspace{1ex}\noindent Proof. Let
$$
v_{\bar\nu}(x_n)=\sum_{j=0}^D(c^0_{\nu,j}+V^0_{\nu,j})x_n^j,\quad
v(x_n)=(v_1(x_n),\dots,v_N(x_n))\;.
$$
It is easy to check that $\delta(x,v(x_n))$ is $x_n$-regular of order $p$.
By the Weierstrass Division Theorem
$$
g_i(x,v(x_n))=Q_i(x,V^0)\delta(x,v(x_n))^2+\sum_{j=0}^{2p-1}R_{i,j}(x',V^0)x_n^j
$$
for $i=1,\dots,m$. Let
$$
\bar v_{\nu}(x)=\sum_{j=0}^D(c^0_{\nu,j}+\bar v^0_{\nu,j})x_n^j,\quad
\bar v(x)=(\bar v_1(x),\dots,\bar v_N(x))
$$
where $(\bar v^0_{\nu,j}(x'))$ is a family of formal power series without
constant term. Thus we get
$$
g_i(x,\bar v(x))=Q_i(x,\bar v(x))\delta(x,\bar v(x))^2+
\sum_{j=0}^{2p-1}R_{i,j}(x',\bar v^0_{\nu,j}(x'))x_n^j\mbox{ for }i=1,\dots,m\;.
$$
By the uniqueness of the remainder in the Weierstrass Division Theorem
we have that $\bar v(x)$ is an approximate solution of the system of
analytic equations $f(x,y)=0$ if and only if $R_{i,j}(x',(\bar v^0_{\nu,j}(x')))=0$
for $i=1,\dots,m$ and $j=0,1,\dots,2p-1$ in $\bbK[[x']]$. This proves
the proposition~\rule{1ex}{1ex}

\section{Proof of the theorem (by induction on the number $n$ of
variables $x$)}

The theorem is trivial for $n=0$. Suppose that $n>0$ and that the theorem
is true for $n-1$. By PROPOSITION~2.1 we may suppose that $\bar y(x)$ is
a simple solution of the system $f(x,y)=0$. Let
$$
\delta(x,y)=\det\frac{J(f_1,\dots,f_m)}{J(y_{N-m+1},\dots,y_N)}\;.
$$
Without diminishing the generality we may suppose that
$\delta(x,\bar y(x))\neq 0$. If $\delta(0,0)\neq 0$ then the theorem follows
from the Implicit Function Theorem. Suppose that $\delta(0,0)=0$.
After a linear change of the variables $x_1,\dots,x_n$ we may assume that
$\delta(x,\bar y(x))$ is $x_n$-regular of order $p>0$. By PROPOSITION~4.1
the system of equations $f(x,y)=0$ has an approximate solution
$\bar v(x)=(\bar v_1(x),\dots,\bar v_N(x))\in\bbK[[x']][x_n]$ such that
the solution $\bar y(x)$ is subordinate to $\bar v(x)$. Write
$$
\bar v_{\nu}(x)=\sum_{j=0}^D(c^0_{\nu,j}+\bar v^0_{\nu,j}(x'))x_n^j,\quad
D\geq 0\mbox{ an integer}
$$
where $(\bar v_{\nu,j}(x'))$ is a family of formal power series without
constant term. It is easy to check that
$$
\left(\sum_{j=0}^Dc_{1,j}^0x_n^j,\dots,
\sum_{j=0}^Dc_{N,j}^0x_n^j\right)
$$
is an approximate solution of the system $f(0,x_n,y)=0$ such that
$$
\ord1\,\delta
\left(0,x_n,\sum_{j=0}^Dc_{1,j}^0x_n^j,\dots,
\sum_{j=0}^Dc_{N,j}^0x_n^j\right)=p\;.
$$
By PROPOSITION~4.2 there exist convergent power series
$F(x',V^0)\in\bbK\{x',V^0\}^M$ such that $F(x',(\bar v^0_{\nu,j}(x')))=0$.
By induction hypothesis there exist convergent power series
$(V^0_{\nu,j}(x',s))$ in $\bbK\{x',s\}$ where $s=(s_1,\dots,s_q)$
are new variables and formal power series
$\bar s(x')=(\bar s_1(x'),\dots,\bar s_q(x'))$ without constant term
such that
$$
F(x',(V^0_{\nu,j}(x',s)))=0,\quad V^0_{\nu,j}(x',\bar s(x'))=
\bar v^0_{\nu,j}(x')\;.
$$
Let
$$
v_{\nu}(x,s)=\sum_{j=0}^D(c^0_{\nu,j}+V^0_{\nu,j}(x',s))x_n^j
\mbox{ for }\nu=1,\dots,N
$$
and $v(x,s)=(v_1(x,s),\dots,v_N(x,s))$. Thus
$\bar v_\nu(x)=v_\nu(x,\bar s(x'))$ for $\nu=1,\dots,N$.
Again by Proposition~4.2 $v(x,s)$ is an approximate solution of the
system $f(x,y)=0$. By the Bourbaki-Tougeron Implicit Function Theorem
the system $f(x,y)=0$ has the parametric solution determined by $v(x,s)$:
\begin{eqnarray*}
y_\nu(x,s,t) & = & v_\nu(x,s)+\delta(x,v(x,s))^2t_\nu\mbox{ for }
\nu=1,\dots,N-m\\
y_\nu(x,s,t) & = & v_\nu(x,s)+\delta(x,v(x,s))u_\nu(x,s,t)\mbox{ for }
\nu=N-m+1,\dots,N\;.
\end{eqnarray*}
On the other hand
\begin{eqnarray*}
\bar y_\nu(x,t) & = & \bar v_\nu(x)+\delta(x,\bar v(x))^2t_\nu\mbox{ for }
\nu=1,\dots,N-m\\
\bar y_\nu(x,t) & = & \bar v_\nu(x)+\delta(x,\bar v(x))\bar u_\nu(x,t)\mbox{ for }
\nu=N-m+1,\dots,N
\end{eqnarray*}
is the parametric solution determined by $\bar v(x)$.
Since the formal solution $\bar y(x)$ is subordinate to the approximate
solution $\bar v(x)$ there exist formal power series
$\bar t(x)=(\bar t_1(x),\dots,\bar t_{N-m}(x))$, $\bar t(0)=0$,
such that $\bar y(x)=\bar y(x,\bar t(x))$. We have
\begin{eqnarray*}
y_\nu(x,,\bar s(x'),t) & = & \bar v_\nu(x)+\delta(x,\bar v(x))^2t_\nu\mbox{ for }
\nu=1,\dots,N-m\\
y_\nu(x,\bar s(x'),t) & = & \bar v_\nu(x)+\delta(x,\bar v(x))u_\nu(x,\bar s(x'),t)\mbox{ for }
\nu=N-m+1,\dots,N
\end{eqnarray*}
By the uniqueness of the parametric solution determined by the approximate
solution $\bar v(x)$ we get
$$
y(x,\bar s(x'),\bar t(x))=\bar y(x,\bar t(x))=\bar y(x)\;~\rule{1ex}{1ex}
$$

\vspace{2ex}\noindent
{\sc Department of Mathematics\\
Kielce University of Technology\\
AL. 1000 L PP 7\\
25-314 Kielce, Poland}

\end{document}